  \documentclass[a4paper,12pt]{amsart} 
  \usepackage{latexsym} 
  \usepackage[all]{xy} 
  \usepackage{amsfonts} 
  \usepackage{amsthm} 
  \usepackage{amsmath} 
  \usepackage{amssymb} 
  \newcommand{\cC}{{\mathcal C}} 

  \def\sw#1{{\sb{(#1)}}} 
  \def\su#1{{\sp{[#1]}}} 

  \def\endproof{\hbox{$\sqcup$}\llap{\hbox{$\sqcap$}}\medskip} 
  \def\<{{\langle}} 
  \def\>{{\rangle}} 
   
  \def\eps{\varepsilon}

  \def\note#1{{}}

  \def\note#1{}

  \def\cM{{\bf M}}

  \def\cC{{\mathcal C}} 
   
  \def\cD{{\mathcal D}} 
   \def\cE{{\mathcal E}}

  \def\lrhom#1#2#3#4{{{\rm Hom}\sb{#1, #2}(#3,#4)}} 
  \def\lhom#1#2#3{{{\rm Hom}\sb{#1-}(#2,#3)}} 
   
\def\hom#1#2#3{{{\rm Hom}\sb{#1}(#2,#3)}} 
  \def\lend#1#2{{{\rm End}\sb{#1-}(#2)}}

  \def\beq{\begin{equation}} 
  \def\eeq{\end{equation}} 
  \def\DC{{\Delta_\cC}} 
  \def \eC{{\eps_\cC}} 
  \def\DD{{\Delta_\cD}} 
  \def \eD{{\eps_\cD}}

  \def\Desc{{\bf Desc}}

  \def\ot{{\otimes}}

  \newcommand{\Ra}{\Rightarrow}

\def\cex{{\bf CrgExt}}

  \newcounter{zlist} 
  \newenvironment{zlist}{\begin{list}{(\arabic{zlist})}{ 
  \usecounter{zlist}\leftmargin2.5em\labelwidth2em\labelsep0.5em 
  \topsep0.6ex
  \parsep0.3ex plus0.2ex minus0.1ex}}{\end{list}}

  \newcounter{blist} 
  \newenvironment{blist}{\begin{list}{(\alph{blist})}{ 
  \usecounter{blist}\leftmargin2.5em\labelwidth2em\labelsep0.5em 
  \topsep0.6ex 
  \parsep0.3ex plus0.2ex minus0.1ex}}{\end{list}} 

  \newcounter{rlist} 
  \newenvironment{rlist}{\begin{list}{(\roman{rlist})}{ 
  \usecounter{rlist}\leftmargin2.5em\labelwidth2em\labelsep0.5em 
  \topsep0.6ex 
  \parsep0.3ex plus0.2ex minus0.1ex}}{\end{list}}

  \def\Label#1{\label{#1}\ifmmode\llap{[#1] }\else 
  \marginpar{\smash{\hbox{\tiny [#1]}}}\fi} 
  \def\Label{\label}

  \newcounter{c} 
   
  \newcommand{\etyk}[1]{\vspace{-7.4mm}$$\begin{equation}\Label{#1} 
  \addtocounter{c}{1}} 
  \renewcommand{\]}{\ifnum \value{c}=1 $$\else \end{equation}\fi} 
  \newtheorem{proposition}{Proposition}[section] 
  \newtheorem{lemma}[proposition]{Lemma} 
   
  \newtheorem{corollary}[proposition]{Corollary} 
  \newtheorem{theorem}[proposition]{Theorem} 

  \theoremstyle{definition} 
  \newtheorem{definition}[proposition]{Definition} 
   
   \newtheorem{examples}[proposition]{Examples}

  \theoremstyle{remark}

\begin{document} 
~\vspace{.8in}
  \title{A note on coring extensions} 
  \author{Tomasz Brzezi\'nski} 
  \address{ Department of Mathematics, University of Wales Swansea, 
  Singleton Park, \newline\indent  Swansea SA2 8PP, U.K.} 
  \email{T.Brzezinski@swansea.ac.uk} 
  \urladdr{http//www-maths.swan.ac.uk/staff/tb} 
  \begin{abstract} 
  A notion of a coring extension is defined and it is related to the existence of an additive functor between comodule categories that factorises through forgetful functors. This correspondence between coring extensions and factorisable functors is illustrated by functors between categories of descent data. A category in which objects are corings and morphisms are coring extensions is also introduced.
    \end{abstract} 
  \maketitle 

  \section{Introduction} 
 Given two algebras $A$, $B$ over a commutative ring $k$, an algebra extension or an algebra map $B\to A$ can be equivalently characterised as a $k$-additive functor $F: \cM_A\to\cM_B$ with the factorisation property
$$ \xymatrix{\cM_A\ar[rd]_{U_A}\ar[rr]^F& & \cM_B \ar[ld]^{U_B}\\
& \cM_k, &}$$
where  $U_A$, $U_B$ are forgetful functors (cf.\ \cite{Par:Ver}). Through this correspondence, morphisms of $k$-algebras can be defined as functors having such a factorisation property. This point of view is taken up in a recent paper by Pareigis \cite{Par:ten}, in which functors between categories of entwined modules are studied, conditions for the factorisation property are derived and these are then suggested as the definition of morphisms between entwining structures. The resulting notion of morphisms of entwining structures is different from the one introduced earlier in \cite{Brz:mod}. Since any entwining structure gives rise to a coring such that the entwined modules can be identified with its   comodules (cf.\ \cite{Brz:str}), it is natural to look at the results of \cite{Par:ten} from the coring point of view. This is the aim of the present note in which, rather than changing the established notion of a morphism of corings (cf.\ \cite{Gom:sep}, \cite[Section~24]{BrzWis:cor}), we introduce the notion of an {\em extension of corings} or a {\em coring 
extension}
and show that such extensions arise from and -- provided they satisfy a suitable purity condition (for example in the case of corings associated to entwining structures) -- give rise to $k$-additive functors with a (suitable) factorisation property.

We work over a commutative associative ring $k$ with a unit. All algebras are over $k$, associative and with a unit. The symbol $\otimes$ between $k$-modules and $k$-module maps means tensor product over $k$. As a rule we do not decorate $\otimes$ between elements, unless there is a danger of confusion.
 For a $k$-algebra $A$, the category of right $A$-modules and right $A$-linear maps is denoted by $\cM_A$. The product map in $A$ is denoted by $\mu_A:A\ot A\to A$ and the unit (either as an element of $A$ or as a $k$-linear map $k\to A$) is denoted by $1_A$.  Given  a $k$-algebra $A$, coproduct in an $A$-coring $\cC$ is denoted by  $\DC :\cC\to \cC\ot_A\cC$, and  the counit is denoted by $\eC:\cC\to A$. We use the Sweedler sigma notation,
i.e., for all $c\in \cC$,
$$
\DC(c) = \sum c\sw1\ot c\sw 2, \qquad (\DC\ot_A\cC)\circ\DC(c) =
(\cC\ot_A\DC)\circ\DC(c) = \sum c\sw1\ot c\sw 2\ot c\sw 3,
$$
etc. For an $A$-coring $\cC$, the {\em left dual ring} is defined as a $k$-module ${}^*\cC = \lhom A \cC A$ with the unit $\eC$ and the product, for all $f,g\in {}^*\cC$, $c\in \cC$,  $f*g(c) = \sum g(c\sw 1f(c\sw 2))$.

The category of right $\cC$-comodules and right $\cC$-colinear maps is denoted by $\cM^\cC$.  For a right $\cC$-comodule $M$, $\varrho^M:M\to M\ot_A\cC$ denotes a coaction.  Recall that $\cM^\cC$ is built upon the category of right $A$-modules, in the sense that every right $\cC$-comodule is a right $A$-module, coactions and morphisms are right $A$-linear maps (with additional compatibility conditions). On elements, $\varrho^M$ is denoted by the Sweedler notation $\varrho^M(m) = \sum m\sw 0\ot m\sw 1$ (but see an exception in the proof of Theorem~\ref{thm.main}).  Similar notational conventions apply to coalgebras and their comodules. A detailed account of the theory of corings and comodules can be found in \cite{BrzWis:cor}.

\section{Extensions of corings and factorisable functors}
\label{sec.main}

Recall that, given an $A$-coring $\cC$ and a $B$-coring $\cD$, a $(\cC,\cD)$-bicomodule is a left $\cC$-comodule that is at the same time a right $\cD$-comodule with $\cC$-colinear  $\cD$-coaction. The $\cC$-colinearity of $\cD$-coaction is equivalent to $\cD$-colinearity of $\cC$-coaction.
\begin{definition}\label{def.exten}
Let $A$ and $B$ be $k$-algebras. A $B$-coring $\cD$ is called a {\em right extension} of an $A$-coring $\cC$ provided $\cC$ is a $(\cC,\cD)$-bicomodule with the left regular coaction $\DC$.
\end{definition}

For example, if $\cC$ and $\cD$ are $A$-corings and $\gamma:\cC\to\cD$ is an $A$-coring morphism, then $\cC$ is a $(\cC,\cD)$-bicomodule with the left regular coaction $\DC$ and the right coaction $\varrho^\cC=(\cC\ot_A\gamma)\circ\DC$. Thus any $A$-coring morphism gives rise to a coring extension.

Definition~\ref{def.exten} implies in particular that if $\cD$ is a right extension of $\cC$, then necessarily $\cC$ is a right $B$-module and $\DC$ is a right $B$-linear map. This leads to the following

\begin{definition}\label{def.mea}
Let $A$ and $B$ be $k$-algebras. An $A$-coring $\cC$ is said to {\em measure $B$ to $A$} if there exists a left $A$-linear map $\nu:\cC\ot B\to A$ rendering commutative the following diagrams:
\begin{blist}
\item $ \xymatrix{\cC\ar[rd]_{\eC}\ar[rr]^{\cC\ot  1_B}& & \cC\ot B \ar[ld]^{\nu}\\
& A, &}$
\item $ \xymatrix{\cC\ot B\ot B\ar[d]_{\DC\ot B\ot B}\ar[rr]^{\cC\ot  \mu_B}&& \cC\ot  B\ar[rr]^\nu&& A \\
\cC\ot_A\cC\ot  B\ot  B\ar[rr]^{\cC\ot_A\nu\ot  B}&& \cC\ot_AA\ot B\ar[rr]^\simeq &&\cC\ot  B\ar[u]_{\nu}.}$
\end{blist}
The map $\nu$ is called a {\em $\cC$-measuring} of $B$ to $A$.
\end{definition}

\begin{proposition}\label{prop.equiv}
Let $\cC$ be an $A$-coring and $B$ an algebra.  $\cC$-measurings of $B$ to $A$ are in bijective correspondence with algebra maps $B\to {}^*\cC$. 
\end{proposition}
\begin{proof} There is a bijective correspondence between left $A$-module maps $\nu: \cC\ot B\to A$ and $k$-linear maps $\chi:B\to {}^*\cC$ provided by the hom-tensor isomorphism
$$\lhom A {\cC\ot B} A\simeq \hom k B{ \lhom A\cC A} = \hom k B {{}^*\cC}.$$
Explicitly, for all $b\in B$, $c\in \cC$, $\nu(c\ot b)= \chi(b)(c)$. 
Since the counit $\eC$ of $\cC$ is the unit in ${}^*\cC$, the map $\chi$ is 
unital if and only if $\nu(c\ot 1_B) = \chi(1_B)(c) = \eC(c)$. Thus the 
unitality of $\chi$ is equivalent to the commutativity of the diagram (a) 
in Definition~\ref{def.mea} for $\nu$. Second, the multiplicativity of $\chi$ 
means that, for all $b,b'\in B$, $\chi(bb') = \chi(b)*\chi(b')$, i.e., for 
all $c\in \cC$, 
$$
\chi(bb')(c) = \sum \chi(b')(c\sw 1\chi(b)(c\sw 2)).
$$
Therefore, $\chi$ is a multiplicative map if and only if
\begin{eqnarray*}
\nu (c\ot bb') &=& \chi(bb')(c) = \sum \chi(b')(c\sw 1\chi(b)(c\sw 2))\\
& =&  \sum \nu(c\sw 1\chi(b)(c\sw 2)\ot b') = 
\sum \nu (c\sw 1 \nu(c\sw 2 \ot b)\ot b'),
\end{eqnarray*}
i.e., the diagram (b) in Definition~\ref{def.mea} is commutative.
\end{proof}

In  view of  Proposition~\ref{prop.equiv},  any {\em $B^{op}|A$-coring}
in the sense of \cite[Definition~3.5]{Tak:mor}, any {\em right rational pairing} 
of corings in the sense of \cite{ElKGom:sem} or a
{\em measuring left $A$-pairing} of \cite{Abu:rat} are examples
of a $\cC$-measuring.
We illustrate the notion of a $\cC$-measuring
with a number of additional examples.

\begin{examples}\label{ex.measure} ~\\\vspace{-.7\baselineskip}
\begin{zlist}
\item An algebra $A$, viewed as a trivial $A$-coring,  
measures $B$ to $A$ if and only 
if there is an algebra map $B\to A$.
\item Given algebras $A$ and $B$, let $\Sigma$ be a $(B,A)$-bimodule 
that is finitely generated and projective as a right $A$-module and let $\cC =\Sigma^*\ot_B\Sigma$ be the corresponding comatrix $A$-coring 
(cf.\ \cite{KaoGom:com}). 
Then $\cC$-measurings of $B$ to $A$ are in bijective correspondence 
with right $B$-module structures on $\Sigma$ that make $\Sigma$ 
a $(B,B)$-bimodule.
\item Given an algebra map $\iota: B\to A$, take $\cC =A\ot_B A$ 
the canonical Sweedler coring. Fix a left $B$-module structure on 
$A$ provided by the map $\iota$, i.e., $ba := \iota(b)a$.
Then $\cC$-measurings of $B$ to $A$ are in bijective correspondence 
with right $B$-module structures on $A$ that make $A$ a 
$(B,B)$-bimodule.
\item Let $A$ be a $k$-algebra and $C$ be a $k$-coalgebra with 
coproduct $\Delta_C$ and counit $\eps_C$. If $\cC=A\ot  C$ is the 
coring associated to an entwining structure $(A,C,\psi)$ then 
$\cC$-measurings  of $B$ to $A$ are in bijective correspondence with 
{\em entwined measurings} in the sense of \cite[Remark~2.2]{Par:ten}, i.e., 
with $k$-linear maps $f: C\ot  B\to A$ making the following diagrams:
$$ 
\xymatrix{C\ar[d]_{\eps_C}\ar[rr]^{C\ot  1_B}& & C\ot B \ar[d]^{f}\\
k\ar[rr]^{1_A}& &A}$$
and
$$ \xymatrix{C\ot B\ot B\ar[d]_{\Delta_C\ot B\ot B}\ar[rrr]^{C\ot  \mu_B}&&& C\ot  B\ar[rrr]^f&&& A \\
C\ot C\ot  B\ot  B\ar[rr]^{C\ot f\ot  B}&& C\ot A\ot B\ar[rr]^{\psi\ot B} &&A\ot  C\ot  B\ar[rr]^{A\ot f}&&A\ot A \ar[u]_{\mu_A}}$$
commute. 
\end{zlist}
\end{examples}
{\sl Check.}
(1)  This follows immediately from Proposition~\ref{prop.equiv}, since ${}^*A \simeq A$ as $k$-algebras.

(2) Recall from \cite{KaoGom:com} that $\cC = \Sigma^*\ot_B\Sigma$ is an $A$-coring with coproduct and counit, for all $s\in \Sigma$, $s^*\in \Sigma^*$,
$$
\DC(s^*\ot_B s) = \sum_{i\in I} s^*\ot_B 
e_i\ot_A e_i^*\ot_B s, \qquad \eC(s^*\ot_B s) = s^*(s),
$$
 where $\{e_i\in \Sigma, e^*_i\in \Sigma^*\}_{i\in I}$ 
 is a finite dual basis of $\Sigma_A$. 
 Recall also that ${}^*\cC \simeq \lend B\Sigma^{op}$, where the product in $\lend B\Sigma$ is given by $fg(s) = f(g(s))$. Therefore, by Proposition~\ref{prop.equiv}, $\cC$-measurings are in bijective correspondence with anti-algebra maps $B\to \lend B\Sigma$, i.e., with right $B$-module structures on $\Sigma$ such that $\Sigma$ is a $(B,B)$-bimodule. 

(3) This is a special case of (2), simply take $\Sigma = A$ and view $A$ as a left $B$-module via the map $\iota$.

(4) Recall that an entwining structure consists of an algebra $(A,\mu_A, 1_A)$, a coalgebra $(C,\Delta_C,\eps_C)$ and a $k$-linear map $\psi: C\ot A\to A\ot C$ satisfying a number of conditions (cf., e.g., \cite[Section~32]{BrzWis:cor}). In this case, $\cC = A\ot C$ is an $A$-bimodule via $a(a'\ot c)a'' = aa'\psi(c\ot a'')$ and it has a coproduct and counit $\DC(a\ot c) = a\ot \Delta_C(c)$, $\eC(a\ot c) =a\eps_C(c)$. In the view of the isomorphism $\lhom A {A\ot C\ot B} A \simeq \hom k{C\ot B} A $ any $\cC$-measuring $\nu: \cC\ot B\to A$ corresponds to a $k$-linear map $f:C\ot B\to A$, via $\nu(a\ot c\ot b) = af(c\ot b)$. With the help of this identification one immediately checks that the diagram (a) in Definition~\ref{def.mea} for $\nu$ is equivalent to the first of diagrams in (3) for $f$.  As to the second pair of diagrams, introduce the explicit notation $\psi(c\ot a) = \sum_\alpha a_\alpha\ot c^\alpha$, take any $c\in C$, $b,b'\in B$, use the diagram (b) in Definition~\ref{def.mea} and the definition of the right $A$-multiplication on $\cC$ to compute
\begin{eqnarray*}
f(c\ot bb') &=& \nu(1_A\ot c\ot bb') = \sum \nu((1_A\ot c\sw 1)\nu(1_A\ot c\sw 2\ot b)\ot b')\\
&=& \sum_\alpha \nu (f(c\sw 2\ot b)_\alpha\ot c\sw 1^\alpha\ot b') = \sum_\alpha f(c\sw 2\ot b)_\alpha f(c\sw 1^\alpha\ot b').
\end{eqnarray*}
This is exactly the contents of the second of the diagrams in (3). Similarly one proves that if $f$ makes this diagram commutative, then also $\nu$ renders commutative the diagram (b) in Definition~\ref{def.mea}.
\endproof

In particular, if in Example~\ref{ex.measure}~(4) the trivial entwining $\psi: C\ot A\to A\ot C$, $c\ot a\mapsto a\ot c$ is taken, then $\cC = A\ot C$ measures $B$ to $A$ if and only if $C$ measures $B$ to $A$ in the sense of Sweedler \cite[p.\ 138]{Swe:Hop} (this justifies the choice of the name). Also, the combination of Example~\ref{ex.measure}~(4) and Proposition~\ref{prop.equiv} leads to an equivalent description of entwined measurings as algebra maps $B\to \#_\psi(C,A)$, where $\#_\psi(C,A)$ is a $\psi$-twisted convolution algebra defined as a $k$-module $\hom k C A$ with the unit $\eps_C$ and with the product, for all $f,g\in \hom k C A$ and $c\in C$,
$$ 
(f \#_\psi g)(c) = \sum_\alpha f(c\sw 2)_\alpha g(c\sw 1^\alpha).
$$
This follows immediately from the fact that $\#_\psi(C,A)$ is isomorphic to the left dual ring of the coring $\cC = A\ot C$ associated to an entwining structure $(A,C,\psi)$.

Examples~\ref{ex.measure} indicate that the notion of a coring measuring can be understood as one that unifies the notions of an algebra map, a bimodule structure  and entwined measuring. The relationship between a measuring  and a coring extension is revealed in the following
\begin{lemma}\label{lemma.mea}
Let $A$, $B$ be $k$-algebras and let $\cC$ be an $A$-coring. Then the following statements are equivalent:
\begin{zlist}
\item there exists a right $B$-module structure on $\cC$ such that $\cC$ is an $(A,B)$-bimodule and the coproduct $\DC$ is a $B$-linear map;
\item  $\cC$ measures $B$ to $A$.
\end{zlist}
\end{lemma}
\begin{proof}
(1) $\Ra$ (2) Suppose that $\cC$ is an $(A,B)$-bimodule with a $B$-linear coproduct $\DC$, and let $\varrho_\cC : \cC\ot  B\to \cC$ be the right $B$-multiplication. Define $\nu = \eC\circ\varrho_\cC : \cC\ot  B\to A$. The condition (a) 
in Definition~\ref{def.mea} 
for $\nu$ follows then from the commutative diagram
$$
 \xymatrix{\cC\ar[d]_\eC\ar@{=}[rrd]\ar[rr]^{\cC\ot  1_B}& & \cC\ot B \ar[d]^{\varrho_\cC}\\
A& & \cC\ar[ll]^\eC ,}
$$
in which the right upper triangle is commutative by the unitality of the multiplication $\varrho_\cC$. Since $\DC$ is a right $B$-module morphism,
$$
 \xymatrix{\cC\ot  B\ar[d]_{\DC\ot  B}\ar[rr]^{\varrho_\cC} & &\cC\ar[d]^{\DC}\ar[rrd]^{\simeq} &&\\
\cC\ot_A\cC\ot  B\ar[rr]^{\cC\ot_A\varrho_\cC}&&\cC\ot_A\cC\ar[rr]^{\cC\ot_A\eC} && \cC\ot_A A}
$$
is a commutative diagram. The right-hand triangle is simply the counit axiom. Note that the composition of maps in the bottow row equals $\cC\ot_A\nu$. Since $\varrho_\cC$ is an associative multiplication and the above diagram commutes, we  obtain the following commutative diagram
$$
 \xymatrix{\cC\ot_A\cC\ot  B\ot  B\ar[d]_{\cC\ot_A\nu\ot  B} && \cC\ot B\ot B\ar[d]_{\varrho_\cC\ot B}\ar[ll]_{\DC\ot B\ot B}\ar[rr]^{\cC\ot  \mu_B}&& \cC\ot  B\ar[rr]^{\varrho_\cC}&& \cC\ar@{=}[dll]\ar[d]^{\eC} \\
\cC\ot_AA\ot B\ar[rr]^\simeq &&\cC\ot  B\ar[rr]^{\varrho_\cC}&&\cC\ar[rr]^{\eC}&& A.}
$$
The outer rectangle in the above diagram is equivalent to the condition (b) in Definition~\ref{def.mea} for $\nu= \eC\circ\varrho_\cC$. Thus we conclude that $\nu$ is a $\cC$-measuring of $B$ to $A$ as required.

(2) $\Ra$ (1) Given a $\cC$-measuring $\nu$, define $\varrho_\cC: \cC\ot  B\to \cC$ by $c\ot b\mapsto \sum c\sw 1\nu(c\sw 2\ot b)$. Then condition (a) in Definition~\ref{def.mea} for $\nu$ implies that, for all $c\in \cC$,
$$
\varrho_\cC(c\ot 1_B) = \sum c\sw 1\nu(c\sw 2\ot 1_B) = \sum c\sw 1\eC(c\sw 2) = c.
$$
Furthermore, the use of condition (b) (in the second equality below),  
the definition of $\varrho_\cC$ in terms of $\nu$, and the right $A$-linearity of
$\DC$ 
give, for all $c\in\cC$, $b,b'\in B$,
\begin{eqnarray*}
\varrho_\cC(c\ot bb') &=& \sum c\sw 1\nu(c\sw 2\ot bb') = \sum c\sw 1\nu(c\sw 2\nu(c\sw 3\ot b)\ot b')\\
& = & \sum \varrho_\cC(c\sw 1\nu(c\sw 2\ot b)\ot b') =  \varrho_\cC(\varrho_\cC(c\ot b)\ot b'). 
\end{eqnarray*} 
Thus $\cC$ is a right $B$-module with the multiplication $\varrho_\cC$. The map $\varrho_\cC$ is a composition of left $A$-linear maps, hence a left $A$-linear map, i.e., $\cC$ is an $(A,B)$-bimodule. Finally, since the coproduct is a coassociative right $A$-linear map, for all $b\in B$ and $c\in\cC$,
$$
\DC(\varrho_\cC(c\ot b)) = \sum\DC(c\sw 1\nu(c\sw 2\ot b)) = \sum c\sw 1\ot c\sw 2\nu(c\sw 3\ot b) = \sum c\sw 1\ot \varrho_\cC(c\sw 2\ot b).
$$
This means that the coproduct $\DC$ is a right $B$-linear map as required.
\end{proof}

 The main result of this 
note is contained in the following
\begin{theorem}\label{thm.main}
Let $\cC$ be an $A$-coring and $\cD$ be a $B$-coring. 
\begin{zlist}
\item If there exists a $k$-additive functor $F:\cM^\cC\to \cM^\cD$ with a factorisation property
$$ \xymatrix{\cM^\cC\ar[rd]_{U^\cC}\ar[rr]^F& & \cM^\cD \ar[ld]^{U^\cD}\\
& \cM_k, &}$$
where  $U^\cC$, $U^\cD$ are forgetful functors, then $\cD$ is a right extension of $\cC$.
\item Let $\cD$ be a right extension of $\cC$ such that, for all $\cC$-comodules $(N,\varrho^N)$, the right $B$-module map $\varrho^N \ot_A \cC - N\ot_A\DC$ is $\cD\ot_B \cD$-pure. Then there exists a $k$-additive functor $F:\cM^\cC\to \cM^\cD$ with a factorisation property as in (1).
\end{zlist}
\end{theorem}

\begin{proof}
(1)  Let $F:\cM^\cC\to \cM^\cD$ be a $k$-additive functor that factorises through forgetful functors $U^\cC$ and $U^\cD$. The factorisation property means that for any $M\in \cM^\cC$, $F(M) = M$ as $k$-modules. Similarly, for any morphism $f$ in $\cM^\cC$, $F(f) =f$ as $k$-linear maps. This implies that for any right $\cC$-comodule $M$ there exist an action $\varrho_M :M\ot  B\to M$ and a coaction $\varrho^M: M\to M\ot_B\cD$, and any $k$-linear map $f$ that is a morphism in $\cM^\cC$ is also a morphism in $\cM^\cD$ (the functoriality of action and coaction). In particular, $\cC$ is a right $\cC$-comodule with the regular coaction $\DC$, hence $\cC$ is a right $B$-module and there exists a right $\cD$-coaction on $\cC$,
$\varrho^\cC:\cC\to \cC\ot_B\cD$. In addition, $\DC$ is a left $A$-module map. Equivalently, for any $a\in A$, the right $\cC$-colinear map $\ell_a:\cC\to \cC$, $c\mapsto ac$ is a morphism in $\cM^\cC$. Thus the $k$-linear map $\ell_a$ is a morphism in $\cM^\cD$, i.e., $\varrho^\cC$ is a left $A$-module map. Furthermore,  $\cC\ot_A\cC$ is a right $\cC$-comodule with the  coaction $\cC\ot_A\DC$, hence it is a right $B$-module and there exists a right $\cD$-coaction,
$\varrho^{\cC\ot_A\cC}:\cC\ot_A\cC\to \cC\ot_A\cC\ot_B\cD$. For any $c\in \cC$, consider a right $\cC$-comodule map $\ell^c: \cC\to \cC\ot_A\cC$, $c'\mapsto c\ot c'$. Since $\DC$ is a right $\cC$-comodule map too, the functoriality of $\cD$-coactions  implies that, for all $c\in \cC$,
$$
(\ell^c\ot_B\cD)\circ\varrho^\cC = \varrho^{\cC\ot_A\cC}\circ\ell^c, \qquad \varrho^{\cC\ot_A\cC}\circ\DC = (\DC\ot_B\cD)\circ\varrho^\cC,
$$
in $\cM_k$. Putting these two together we obtain,
\begin{eqnarray*}
\sum c\sw 1\ot \varrho^\cC(c\sw 2) &=& \sum (\ell^{c\sw 1}\ot_B\cD)\circ \varrho^\cC(c\sw 2) = \sum  \varrho^{\cC\ot_A\cC}\circ\ell^{c\sw 1}(c\sw 2)\\
&=& \sum \varrho^{\cC\ot_A\cC}(c\sw 1 \ot c\sw 2) = \varrho^{\cC\ot_A\cC}\circ\DC(c)\\
&=& (\DC\ot_B\cD)\circ\varrho^\cC(c).
\end{eqnarray*}
This means that the coaction $\varrho^\cC$ is a left $\cC$-comodule map, hence $\cC$ is a $(\cC,\cD)$-bicomodule, i.e., $\cD$ is a coring extension of $\cC$.

(2) This is contained in \cite[22.3, Erratum]{BrzWis:cor}. In detail, suppose that $\cD$ is a coring extension of $\cC$ and let $\nu : \cC\ot B\to A$ be the measuring corresponding to the right $B$-multiplication on $\cC$ as in  Lemma~\ref{lemma.mea}. Write $\sigma: \cC\to \cC\ot_B\cD$ for the right $\cD$-coaction on $\cC$. Define a $k$-linear functor $F:\cM^\cC\to\cM^\cD$ as follows.

Take any right $\cC$-comodule $M$ and define a map $\varrho_M: M\ot B\to M$, $m\ot b\mapsto \sum m\sw 0\nu(m\sw 1\ot b)$. Following the same steps as in the proof (2) $\Ra$ (1) of Lemma~\ref{lemma.mea} one easily verifies that $M$ is a right $B$-module with multiplication $\varrho_M$. We write $m. b := \varrho_M(m\ot b)$. Furthermore, if $f:M\to N$ is a morphism in $\cM^\cC$, then for all $m\in M$, $b\in B$,
\begin{eqnarray*}
f(m.b) &=& \sum f(m\sw 0\nu(m\sw 1\ot b)) = \sum f(m\sw 0)\nu(m\sw 1\ot b)\\
&=& \sum f(m)\sw 0\nu(f(m)\sw 1\ot b) = f(m).b.
\end{eqnarray*} 
The second equality follows from the $A$-linearity of $f$, while the third one is a consequence of the fact that $f$ is a $\cC$-comodule map. The first and last equalities follow from the definition of the $B$-multiplication $\varrho_M$ on $M$. Therefore, $f$ is a right $B$-linear map, and thus we have constructed a functor $\cM^\cC\to \cM_B$. Now we need to define a $\cD$-coaction on any $\cC$-comodule.

Start with the right $\cC$-comodule isomorphism $M\simeq M\Box_\cC \cC$ (cf.\ \cite[22.4]{BrzWis:cor}) provided by the right $\cC$-coaction on $M$. Here $\Box_\cC$ denotes the cotensor product over $\cC$. By applying the above functor $\cM^\cC\to \cM_B$ we obtain a right $B$-module isomorphism and we can construct a map
$$
\varrho^M : \xymatrix@R=2pt{M\ar[r]^(.4){\simeq}&  M\Box_\cC\cC \ar[rr]^(.35){M\Box_\cC\sigma} &&M\Box_\cC(\cC\ot_B\cD)\simeq M\ot_B\cD\\
m\ar@{|->}[rrr]&&&\sum m\sw 0\eC(m\sw 1\su 0)\ot m\sw 1\su 1,}
$$
where $\sum m\sw 0\ot m\sw 1\in M\ot_A\cC$ denotes the $\cC$-coaction on $M$, while $\sigma(c) = \sum c\su 0\ot c\su 1\in \cC\ot_B\cD$ denotes the $\cD$-coaction on $\cC$. 
Note that $\varrho^M$ is well-defined by the purity assumption.
We claim that $\varrho^M$ is a right $\cD$-coaction. First, $\varrho^M$ is a right $B$-module map as a composition of $B$-module maps. Note that, for all $m\in M$,
$$
\sum m\sw 0\eC(m\sw 1\su 0)\eD(m\sw 1\su 1) = \sum m\sw 0\eC(m\sw 1) = m,
$$
so that $\varrho^M$ is a counital map. It remains to check the coassociativity of $\varrho^M$. This is done by a rather lengthy but straightforward calculation, the details of which are displayed below. Take any $m\in M$ and compute
\begin{eqnarray*}
(\varrho^M&&\!\!\!\!\!\!\!\!\!\!\!\!\!\!\ot_B\cD)\circ\varrho^M(m) = \sum \varrho^M (m\sw 0\eC(m\sw 1\su 0))\ot m\sw 1\su 1 \\
&=& \sum m\sw 0\sw 0\eC((m\sw 0\sw 1\eC(m\sw 1\su 0))\su 0)\ot (m\sw 0\sw 1\eC(m\sw 1\su 0))\su 1\ot m\sw 1\su 1\\
&=& \sum m\sw 0\eC((m\sw 1\sw 1\eC(m\sw 1\sw 2\su 0))\su 0)\ot (m\sw 1\sw 1\eC(m\sw 1\sw 2\su 0))\su 1\ot m\sw 1\sw 2\su 1\\
&=& \sum m\sw 0\eC((m\sw 1\su 0\sw 1\eC(m\sw 1\su 0\sw 2))\su 0)\ot (m\sw 1\su 0\sw 1\eC(m\sw 1\su 0\sw 2))\su 1\ot m\sw 1\su 1\\
&=& \sum m\sw 0\eC(m\sw 1\su 0\su 0)\ot m\sw 1\su 0\su 1\ot m\sw 1\su 1\\
&=& \sum m\sw 0\eC(m\sw 1\su 0)\ot m\sw 1\su 1\sw 1\ot m\sw 1\su 1\sw 2\\
&=& (\cC\ot_B\DD)\circ\varrho^M(m).
\end{eqnarray*}
The fourth equality is a consequence of the fact that $\cC$ is a $(\cC,\cD)$-bicomodule
and the purity assumption is used to derive the penultimate equality.
Thus $\varrho^M$ is a right $\cD$-coaction.

We already know that if $f:M\to N$ is a morphism in $\cM^\cC$, then $f$ is a right $B$-module map. Take any $m\in M$ and compute
\begin{eqnarray*}
\varrho^N(f(m)) &=& \sum f(m)\sw 0\eC(f(m)\sw 1\su 0)\ot f(m)\sw 1\su 1 = \sum f (m\sw 0)\eC(m\sw 1\su 0)\ot m\sw 1\su 1\\
&=& \sum f (m\sw 0\eC(m\sw 1\su 0))\ot m\sw 1\su 1 = (f\ot_B\cD)\circ\varrho^M(m),
\end{eqnarray*}
where the second equality follows from the $\cC$-colinearity and the third one from the $A$-linearity of $f$. Thus $f$ is a morphism of right $\cD$-comodules. Put together, all this means that we have constructed a $k$-additive functor $F:\cM^\cC\to \cM^\cD$. It is obvious from this construction  that $F$ factorises through the forgetful functors $U^\cC$ and $U^\cD$ as required.
\end{proof}

In  view of Lemma~\ref{lemma.mea}, a
trivial $B$-coring $B$ is a right 
extension of $\cC$ if and only if $\cC$ measures $B$ to $A$. Thus Theorem~\ref{thm.main} leads immediately to the following characterisation
of measurings in terms of functors with a factorisation property.
\begin{corollary}
Let $A$ and $B$ be $k$-algebras and let $\cC$ be an $A$-coring. The following statements are equivalent:
\begin{zlist}
\item $\cC$ measures $B$ to $A$;
\item there exists a $k$-additive functor $F:\cM^\cC\to \cM_B$ with a factorisation property
$$ \xymatrix{\cM^\cC\ar[rd]_{U^\cC}\ar[rr]^F& & \cM_B \ar[ld]^{U_B}\\
& \cM_k, &}$$
where  $U^\cC$, $U_B$ are forgetful functors.
\end{zlist}
\end{corollary}

If $\cC$ and $\cD$ are corings corresponding to entwining structures $(A,C,\psi)$, $(B,D,\Psi)$ then Theorem~\ref{thm.main} implies \cite[Theorem~2.3]{Par:ten}. As many general results about corings, Theorem~\ref{thm.main} finds an application in non-commutative descent theory.

Given an algebra extension $\iota: B\to A$, the category of right comodules of the associated Sweedler coring $\cC = A\ot_B A$ is isomorphic to the category of {\em (right) descent data} $\Desc (A| B)$ (cf.\ \cite[25.4]{BrzWis:cor}). The objects in $\Desc (A|B)$ are pairs $(M,f)$, where $M$ is a right $A$-module 
and $f: M\to M\ot_B A$ is a right $A$-module map rendering commutative the following diagrams
$$
\xymatrix{M\ar[r]^{f}\ar@{=}[dr] & M\ot_B A \ar[d]^{\varrho_{M|B}}\\
 &M,} \quad  
\xymatrix{M\ar[r]^{f}\ar[dd]_f&  M\ot_B A\ar[rd]^{f\ot_B A} &\\
&& M\ot_B A\ot_B A.\\
M\ot_B A\ar[r]^\simeq& M\ot_BB\ot_B A \ar[ru]_{M\ot_B\iota\ot_B A}&}$$
Here $\varrho_{M|B}: M\ot_B A\to M$ is the 
factorised (through $M\ot A\to M\ot_B A$) $A$-multiplication on $M$. Using this identification of right comodules of  the Sweedler coring $A\ot_BA$ with (right) descent data of the ring extension $B\to A$ we can thus derive the following corollary of Theorem~\ref{thm.main}.
\begin{corollary}\label{cor.descent}
Let $A$, $B$, $D$ be $k$-algebras and let $\iota_B :D\to B$ and $\iota_A: B\to A$ be algebra maps. Then the following statements are equivalent.
\begin{zlist}
\item There exists a $k$-additive functor $F:\Desc (A|B)\to\Desc (B|D)$ with a factorisation property
$$ \xymatrix{\Desc (A|B)\ar[rd]_{U^{A|B}}\ar[rr]^F& & \Desc (B|D) \ar[ld]^{U^{B|D}}\\
& \cM_k, &}$$
where  $U^{A|B}$, $U^{B|D}$ are forgetful functors.
\item View $A$ as a left $B$-module via the map $\iota_A$. There exist:
\begin{rlist}
\item a left $B$-linear right $B$-multiplication $\varrho_A:A\ot B\to A$;
\item a  $(B,B)$-bimodule map $\varphi:A\to A\ot_B A\ot_D B$ ($A$ is a right $B$-module via $\varrho_A$) rendering commutative the following three diagrams:\\
(a) $
\xymatrix{A\ar[rr]^{\simeq}\ar[d]_\varphi & & B\ot_B A \ar[d]^{\iota_A\ot_B A}\\
A\ot_B A\ot_D B\ar[rr]^{A\ot_B\varrho_{A|D}} && A\ot_B A,}$\\
(b) $
\xymatrix{A\ar[r]^{\varphi}\ar[d]_\varphi & A\ot_B A\ot_D B\ar[rr]^{A\ot_B\varphi\ot_DB}& & A\ot_B A\ot_B A\ot_D B\ot_D B \ar[d]_{\mu_{A|B}\ot_B A\ot_D B\ot_D B}\\
A\ot_B A\ot_D B\ar[r]^{\simeq} & A\ot_B A\ot_D D\ot_D B\ar[rr]^{A\ot_B A\ot_D\iota_B\ot_D B} && A\ot_B A\ot_D B\ot_D B,}$\\
(c) $
\xymatrix{A\ar[r]^{\simeq}\ar[d]_\varphi & B\ot_B A \ar[rr]^{\iota_A\ot_B \varphi} && A\ot_BA\ot_B A\ot_D B\\
A\ot_B A\ot_D B\ar[rrr]^{\simeq} &&& A\ot_B B\ot_B A\ot_D B\ar[u]_{A\ot_B\iota_A\ot_BA\ot_D B},}$\\~\\
where $\varrho_{A|D}: A\ot_D B\to A$ is the  factorised (through $A\ot B\to A\ot_D B$) right $B$-multiplication $\varrho_A$ and $\mu_{A|B}: A\ot_B A\to A$ is the  factorised (through $A\ot A\to A\ot_BA$) product $\mu_A$ in $A$. 
\end{rlist}
\end{zlist}
\end{corollary}
\begin{proof}
The purity condition in Theorem~\ref{thm.main}~(2) always holds for  Sweedler corings, hence 
part (1) is equivalent to the statement that $A\ot_B A$ is a right coring extension of $B\ot_D B$. The contents of this statement is contained in part (2). We only explain the origin of maps and diagrams in part (2), the details of the proof are left to the reader. By Example~\ref{ex.measure}(3), there must exist a right $B$-multiplication $\varrho_A$ as in (i). Furthermore, $A\ot_BA$ must be a right $B\ot_DB$-comodule, thus there exists a descent datum $(A\ot_BA,f)\in \Desc (B|D)$. The map $f: A\ot_B A\to A\ot_BA\ot_DB$ is a left $A$-module map as it corresponds to a coaction that is left $\cC$-colinear. As a part of a descent datum, $f$ is right $B$-linear. In view of the isomorphism $\lrhom AB{A\ot_BA}{A\ot_BA\ot_DB}\simeq \lrhom BB{A}{A\ot_BA\ot_DB}$,  $f$ can be equivalently given as a $(B,B)$-bimodule map $\varphi: A\to A\ot_BA\ot_DB$, $a'\varphi(a) = f(a'\ot a)$. The diagrams (a) and (b) are the defining diagrams for $f$ as a part of a descent datum written in terms of $\varphi$. The diagram (c) expresses the fact that $f$ is a left $A\ot_B A$-comodule map, as it corresponds to a coaction that makes $A\ot_B A$ into an $(A\ot_B A,B\ot_DB)$-comodule.
\end{proof}

A functor obtained as a composition of any two functors between comodule categories that factorise through the forgetful functors also factorises through forgetful functors. The correspondence between such functors and coring extensions leads therefore to a category $\cex_k^r$ in which objects are corings understood as pairs $(\cC\!:\! A)$. Morphisms $(\cC\!:\! A)\to (\cD\!:\! B)$ are {\em pure coring extensions}, i.e.\  pairs $(\varrho_{\cC}, \varrho^{\cC})$ where $\varrho_{\cC}: \cC\ot B\to \cC$ is a left $\cC$-colinear $B$-action and $\varrho^{\cC} :\cC\to\cC\ot_B\cD$ is a left  $\cC$-colinear $\cD$-coaction, such that, for all $\cC$-comodules $(N,\varrho^N)$, the right $B$-module map $\varrho^N \ot_A \cC - N\ot_A\DC$ is $\cD\ot_B \cD$-pure. A composition of morphisms $(\cC\!:\! A)\to (\cD\!:\! B)$ and $(\cD\!:\! B)\to (\cE\!:\! D)$ is derived from the composition of corresponding functors and comes out as
$$
\varrho_{\cD}\bullet\varrho_{\cC} : \xymatrix@R=2pt{\cC\ot D\ar[rr]^(.4){\varrho^{\cC}\ot D}& & \cC\ot_B\cD\ot D \ar[rr]^(.55){\cC\ot_B \varrho_{\cD}} &&\cC\ot_B\cD\ar[rr]^{\cC\ot_B\eD}&& \cC\ot_B B \simeq \cC}
$$
and 
$$
\varrho^{\cD}\bullet\varrho^{\cC}  : \xymatrix@R=2pt{\cC\ar[r]^(.4){\simeq}&  \cC\Box_\cD\cD \ar[rr]^(.35){\cC\Box_\cD\varrho^{\cD}} &&\cC\Box_\cD(\cD\ot_D\cE)\simeq \cC\ot_D\cE,}
$$
where the first isomorphism is provided by the coaction $\varrho^{\cC|\cD}$, while the second one is obtained with the help of the counit in $\cD$ (compare the construction of  coaction $\varrho^M$ in the proof of Theorem~\ref{thm.main}~(2)).

Finally, we would like to point out that the results of this note can also be presented for left comodules of a coring thus leading to the notions of a left $\cC$-measuring and a left coring extension. This is achieved by using the obvious left-right correspondence. Note, however, that a right coring extension is not necessarily a left coring extension, thus the left-right symmetry that exists for characterisation of algebra (or coalgebra) extensions does not exist in the general coring case.

  \section*{Acknowledgements} 
  I would like to thank Laiachi El Kaoutit for pointing out to me the relationship
  between measurings and pairings, and 
  the Engineering and Physical 
  Sciences Research 
  Council for an Advanced Fellowship. I am also grateful to Gabriella B\"ohm for spotting a mistake in the original version of this note.

  \end{document}